
\magnification=\magstephalf
\hoffset=39.5 truept\voffset=5 truept
\hsize=30truepc\vsize=53truepc
\parindent=4.3 truemm
\let\sl\it

\font\FSmark=cmr8 scaled 910
\font\FSjournal=cmcsc10 at 6.4pt
\font\FStitle=cmbx10
\font\FSauthor=cmr8
\font\FSthanks=cmr8

\font\FStitleAB=cmcsc10 at 8.22pt
\font\FStextAB=cmr8
\font\FSsubsec=cmcsc10
\font\FSsection=cmcsc10 scaled\magstephalf
\font\FSauthorB=cmcsc10 at 8.22pt
\font\FSboldB=cmbx8
\font\FSromanB=cmr8
\font\FSslantB=cmti8
\font\FSeightrm=cmr8 \font\FSfiverm=cmr5
\font\FSeighti=cmmi8 \font\FSfivei=cmmi5
\font\FSeightsy=cmsy8 \font\FSfivesy=cmsy5
\font\FSeightex=cmex10 at 8pt
\font\eightmsbm=msbm9
\def\ESshortfam{%
\textfont0=\FSeightrm \scriptfont0=\FSfiverm \scriptscriptfont0=\FSfiverm
\textfont1=\FSeighti \scriptfont1=\FSfivei \scriptscriptfont1=\FSfivei
\textfont2=\FSeightsy \scriptfont2=\FSfivesy \scriptscriptfont2=\FSfivesy
\textfont3=\FSeightex \scriptfont3=\FSeightex \scriptscriptfont3=\FSeightex
\textfont\msbmfam=\eightmsbm
\scriptfont\msbmfam=\fivemsbm
\scriptscriptfont\msbmfam=\fivemsbm}%
\newfam\msbmfam 
\font\tenmsbm=msbm10\textfont\msbmfam=\tenmsbm
\font\sevenmsbm=msbm7 \scriptfont\msbmfam=\sevenmsbm
\font\fivemsbm=msbm5 \scriptscriptfont\msbmfam=\fivemsbm 
\mathchardef\N="7\the\msbmfam4E 
\mathchardef\Z="7\the\msbmfam5A 
\mathchardef\Q="7\the\msbmfam51 
\mathchardef\R="7\the\msbmfam52 
\mathchardef\C="7\the\msbmfam43 
\mathchardef\E="7\the\msbmfam45 
\mathchardef\P="7\the\msbmfam50 
\def\ESmarkeven#1{\def\ESmarkE{#1}}%
\def\ESmarkodd#1{\def\ESmarkO{#1}}%
\newdimen\testheadline
\headline={%
\ifnum\pageno=1%
  {\FSjournal ESAIM: Probability and Statistics}%
  \hfill%
\else
  \ifodd\pageno%
\setbox0\hbox{\FSmark\ESmarkO}
\setbox1\hbox{\FSmark\folio}
\testheadline=\wd1
\multiply\testheadline by 2%
\advance\testheadline by \wd0%
\advance\testheadline by 20pt\relax%
\ifdim\testheadline>\hsize
\message{*** WARNING: your short title is too long! ***}\fi
\hfill\box0\hfill\llap{\box1}%
  \else%
\setbox0\hbox{\FSmark\ESmarkE}
\setbox1\hbox{\FSmark\folio}
\testheadline=\wd1
\multiply\testheadline by 2%
\advance\testheadline by \wd0%
\advance\testheadline by 20pt\relax%
\ifdim\testheadline>\hsize
\message{*** WARNING: too many authors in headline! ***}\fi
\rlap{\box1}\hfill\box0\hfill%
  \fi%
\fi}%
\footline={%
\ifnum\pageno=1%
  \hfill%
\else
  \ifodd\pageno%
    \hfill{\FSjournal ESAIM: P\&S}%
  \else%
    {\FSjournal ESAIM: P\&S}\hfill%
  \fi%
\fi}%
\def\EStitle#1{\par\noindent\line{\hfill{\FStitle #1}\hfill}}%
\def\ESauthor#1{\vskip 16pt\noindent\line{\hfill{\FSauthor #1}\hfill}\par}%
\long\def\ESthanks#1{{\FSthanks\baselineskip=10 pt\footnote{}{#1}}}%
\long\def\ESabstract#1{\vskip 16pt\noindent%
\leftskip=37.5 truept\rightskip=37.5 truept\baselineskip=10 pt%
{\FStitleAB #1.}\begingroup\ESshortfam\FStextAB{}}%
\long\def\ESendabstract{\endgroup\par\noindent%
\leftskip=0 truept\rightskip=0 truept\baselineskip=12 pt}%
\newcount\sectionnumber\sectionnumber=0
\newcount\subsectionnumber
\def\ESsection#1{\begingroup%
\global\advance\sectionnumber by 1%
\global\subsectionnumber=0\global\formulanumber=1\global\theoremnumber=1%
\long\def\lf##1\cr{\line{\hfill{\FSsection ##1}\hfill}\par\noindent}%
\goodbreak\vskip 8truemm\noindent\lf\the\sectionnumber. #1\cr
\vskip 2truemm\endgroup\noindent
\def\verif{\ifx\next\par\let\next\ignorepar\fi
\ifx\next\ESsubsection\vskip-2truemm\fi
\next}%
\def\ignorepar{\afterassignment\verif\let\next=}%
\ignorepar}
\def\ESsubsection#1{\begingroup%
\global\advance\subsectionnumber by 1%
\long\def\lf##1\cr{\line{\hfill{\FSsubsec ##1}\hfill}\par\noindent}%
\goodbreak\vskip 2truemm\noindent\lf\the\sectionnumber.\the\subsectionnumber.
#1\cr%
\vskip 2truemm\endgroup\noindent
\def\verif{\ifx\next\par\let\next\ignorepar\fi\next}%
\def\ignorepar{\afterassignment\verif\let\next=}%
\ignorepar}
\newcount\formulanumber
\def\neweqno{\xdef\lasteq{\the\sectionnumber.\the\formulanumber}%
\eqno(\lasteq)\global\advance\formulanumber by 1}
\def\label#1{{\def\next{#1}%
\global\expandafter
\let\csname formula = \meaning\next\endcsname=\lasteq}}
\def\formula#1{{\def\next{#1}%
(\csname formula = \meaning\next\endcsname)}}
\newcount\theoremnumber
\def\newthmno{\xdef\lastthm{\the\sectionnumber.\the\theoremnumber}%
\lastthm\global\advance\theoremnumber by 1}
\def\thmname#1{{\def\next{#1}%
\global\expandafter
\let\csname theorem = \meaning\next\endcsname=\lastthm}}
\def\quotethm#1{{\def\next{#1}%
\csname theorem = \meaning\next\endcsname}}
\def\ESdefinition#1{\vskip 1truemm\noindent{\FSsubsec #1 \newthmno.}}
\def\ESenddefinition{\par}

\def\ESexample#1{\vskip 1truemm\noindent{\FSsubsec #1 \newthmno.}}
\def\ESendexample{\par}
\def\EStheorem#1{\vskip 1truemm\noindent{\FSsubsec #1
\newthmno.}\begingroup\sl}%
\def\ESendtheorem{\endgroup\par\noindent}%
\def\ESproof#1{\smallskip\noindent{\sl #1.}}%
\def\ESendproof{\hfill$\sqcup\mskip-12mu\sqcap$\smallskip\noindent}%
\def\ESacknow#1{\goodbreak\vskip 6truemm%
\centerline{\FSsection #1}
\vskip 2truemm}
\long\def\ESbiblio#1{\par\noindent%
\begingroup\ESacknow{#1}%
\baselineskip=10 pt\leftskip=39 truept\parindent=-39 truept%
\ESshortfam\def\rm{\FSromanB}\def\sc{\FSauthorB}%
\def\sl{\FSslantB}\def\bf{\FSboldB}\rm{}}%
\long\def\ESendbiblio{\par\endgroup\baselineskip=12 pt}%
\long\def\ESend{\vskip 8 pt\noindent
\ESaddresS%
\ifodd\pageno\vfill\eject\fi}%
\def\footnoterule{\kern-3pt\hrule width .8truein\kern 2.6pt}%
\pageno=1\null\vskip 8pt\noindent%
\ESthanks{URL address of the journal: \FSboldB
http://www.emath.fr/ps/}%
\font\mail=cmtt8

\def\dli{\ifhmode\par\fi\noindent}
\def\II{\mathrel 1\joinrel\mathrel {\rm I}}
\def\relmont#1#2{\mathrel{\mathop{\kern 0pt#1}\limits^{\hbox{
$\scriptstyle #2$}}}}


\ESmarkeven{GUY MOREL}

\ESmarkodd{{\it P}-VALUES ET EXPERTISES }
\EStitle{LES {\it P}-VALUES COMME VOTES D'EXPERTS$^{*}$}
\ESauthor{GUY MOREL$^1$}

{\parindent=0mm 
\footnote{}{\parindent=0mm
{\FSslantB Mots cl\'es.} Th\'eorie de la d\'ecision, tests, {\FSslantB p}-values, seuils minimum de rejet, hypoth\`eses unilat\'erales et bilat\'erales.

\FSromanB $^*$ Recherche r\'ealis\'ee dans le cadre du LAST et du CNRS UPRES-A 6083 de Tours.

$^1$ Universit\'e de Tours, UFR Arts et Sciences Humaines, 3 Rue des Tanneurs, 37041 Tours Cedex, France ; e-mail : {\mail morel@univ-tours.fr}
}}

\ESabstract{abstract} The {\it p}-values are often implicitly used as a measure of evidence for the hypotheses of the tests. This practice has been analyzed with different approaches. It is generally accepted for the one-sided hypothesis problem, but it is often criticized for the two-sided hypothesis problem. We analyze this practice with a new approach to statistical inference. First we select good decision rules without using a loss function, we call them experts. Then we define a probability distribution on the space of experts. The measure of evidence for a hypothesis is the inductive probability of experts that decide this hypothesis.
\ESendabstract 

\ESabstract{R\'esum\'e} Dans la pratique des tests, les {\it p}-values sont souvent utilis\'ees comme des mesures du niveau de confiance \`a accorder aux hypoth\`eses. Analys\'ee de diff\'erents points de vue, cette mani\`ere de faire est g\'en\'e\-ra\-lement accept\'ee pour les hypoth\`eses unilat\'erales. Elle est par contre souvent critiqu\'ee dans le cas d'hypoth\`eses bilat\'erales. Nous reprenons ce d\'ebat en utilisant une approche nouvelle des probl\`emes de d\'ecision. Elle consiste \`a s\'electionner de bonnes r\`egles de d\'ecision, appel\'ees experts, sans utiliser de fonction de perte. L'espace des experts est ensuite probabilis\'e. Le poids des experts qui d\'ecident une hypoth\`ese est pris comme indice de confiance en faveur de cette hypoth\`ese.
\ESendabstract 

\ESabstract{AMS Subject Classification.} 62A99, 62C05, 62P.
\ESendabstract

\bigskip
\ESsection{Introduction}%
\medskip

Tous les logiciels statistiques fournissent le r\'esultat d'un test en donnant la {\it p}-value (le niveau de signification). Ce seuil minimum de rejet permet de conclure quel que soit le seuil choisi. Si la suite du programme ne d\'epend pas du r\'esultat du test, l'utilisateur n'est pas oblig\'e de fournir un seuil. On lui \'evite ainsi bien des t\^atonnements, car dans beaucoup d'applications le seuil n'est pas compl\`etement fix\'e. Plusieurs valeurs standards sont envisageables, par exemple $\alpha=0.05$, $\alpha=0.01$ et $\alpha=0.001$.
Lorsque le rejet de l'hypoth\`ese $H_0$ est possible, le r\'esultat est souvent symbolis\'e par des \'etoiles : * signifie que seul $\alpha=0.05$ permet de rejeter $H_0$, *** indique que m\^eme $\alpha=0.001$ le permet. Cette symbolique des \'etoiles n'est bien s\^ur pas innocente, elle traduit que le rejet est d'autant plus s\^ur que les \'etoiles sont nombreuses. Donner la {\it p}-value quand elle est inf\'erieure \`a $0.05$ se fait de plus en plus. Ce n'\'etait pas possible lorsque le r\'esultat d'un test s'obtenait par lecture sur une table, l'utilisateur ne disposait que des valeurs critiques du test pour les seuils tabul\'es. Avec la {\it p}-value la symbolique des \'etoiles se transforme en une phrase du type : plus la {\it p}-value est faible, plus la ``probabilit\'e'' de l'hypoth\`ese $H_0$ est faible. Souvent cette interpr\'etation probabiliste s'\'elargit \`a toute la plage des valeurs possibles, l'intervalle $[0,1]$. La {\it p}-value est alors consid\'er\'ee comme la ``probabilit\'e'' de l'hypoth\`ese $H_0$, son compl\'ement \`a $1$ comme la ``probabilit\'e'' de $H_1$. Les guillemets posent le probl\`eme du statut de cette probabilit\'e, les auteurs de la th\'eorie des tests n'avaient pas pour but de fournir une probabilit\'e pour les deux hypoth\`eses. Les {\it p}-values sont un indice du niveau de confiance accord\'e \`a l'hypoth\`ese $H_0$ qui a \'et\'e forg\'e par la pratique des tests. Cet indice statistique non th\'eoris\'e suscite de nombreux articles. On peut les diviser en deux grandes cat\'egories.

Certains articles partent de l'outil {\it p}-value tel qu'il est construit par la pratique et \'etudient ses propri\'et\'es. Parmi les plus r\'ecents, citons par exemple l'article de Hung, O'Neill, Bauer et K\"ohne (1997) qui traite de la distribution des {\it p}-values sous l'hypoth\`ese $H_1$, et celui de Thompson (1996) qui utilise les {\it p}-values, \`a la place de la puissance, pour comparer les tests.

Les autres articles traitent plut\^ot du fondement des {\it p}-values comme pro\-ba\-bilit\'e sur l'espace des deux hypoth\`eses. Ils les \'etudient dans un autre cadre que la th\'eorie des tests. C'est aussi ce que nous essaierons de faire ici. \dli Le cadre th\'eorique qui para\^{\i}t le plus appropri\'e pour ce type d'\'etude est le cadre bay\'esien : les {\it p}-values peuvent-elles \^etre consid\'er\'ees comme des probabilit\'es a posteriori ? Pour r\'esumer, disons que ce n'est g\'en\'eralement pas le cas, mais que pour les tests unilat\'eraux les {\it p}-values s'obtiennent souvent \`a partir de probabilit\'es a priori non informatives [Berger (1985) ; Berger et Sellke (1987) ; Casella et Berger (1987) ; Robert (1992)].
\dli On peut aussi \'etudier les {\it p}-values comme estimateurs de l'indicatrice de l'ensemble des param\`etres d\'efinissant $H_0$ [Hwang, Casella, Robert, Wells et Farrell (1992) ; van der Meulen et Schaafsma (1993) ; Robert (1992) ; Schaafsma, Tolboom et van der Meulen (1989)]. L\`a encore, pour les tests unilat\'eraux les {\it p}-values sont g\'en\'eralement admissibles, alors qu'elles ne le sont pas pour les tests bilat\'eraux.
\dli Cette opposition entre hypoth\`eses unilat\'erales et hypoth\`eses bilat\'erales se retrouve quand on regarde la coh\'erence des {\it p}-values lorsqu'on fait varier l'hypoth\`ese $H_0$, sans changer le mod\`ele statistique. Consid\'erons deux hypoth\`eses possibles $H_0 : \theta\in\Theta_0$ ($\Theta=\Theta_0+\Theta_1$) et $H'_0 : \theta\in\Theta'_0$ ($\Theta=\Theta'_0+\Theta'_1$). Si $\Theta_0\subset\Theta'_0$ on cherche g\'en\'eralement une mesure du niveau de confiance qui donne une valeur plus grande pour $\Theta'_0$ que pour $\Theta_0$. Ce n'est pas le cas des {\it p}-values pour des hypoth\`eses bilat\'erales, Schervish (1996) le montre dans le mod\`ele statistique d\'efini par la famille des lois gaussiennes $\{N(\theta,1)\}_{\theta\in\R}$. Il trouve facilement $\theta_1<-\theta_0<\theta_0<\theta_2<x$, tels que les hypoth\`eses $[\theta_0]\subset[-\theta_0,\theta_0]\subset[\theta_1,\theta_2]$ et la r\'ealisation $x$ donnent des {\it p}-values qui sont d\'ecroissantes au lieu d'\^etre croissantes. Par exemple, pour $\theta_1=-0.82$, $\theta_0=0.50$, $\theta_2=0.52$ et $x=2.18$ il obtient les {\it p}-values suivantes : $0.0930$ si $\Theta_0=[\theta_0]$, $0.0502$ si $\Theta_0=[-\theta_0,+\theta_0]$ et $0.0498$ si $\Theta_0=[\theta_1,\theta_2]$. Ce sont les valeurs respectives de $2[1-\Phi(2.18-0.50)]$, $[1-\Phi(2.18-0.50)]+[1-\Phi(2.18+0.50)]$ et $[1-\Phi(2.18-0.52)]+[1-\Phi(2.18+0.82)]$, $\Phi$ \'etant la fonction de r\'epartition de la loi normale r\'eduite (voir le paragraphe 4.2).

Les {\it p}-values peuvent aussi \^etre retrouv\'ees \`a partir de la distribution fiduciaire de Fisher mais ce dernier tenait \`a bien distinguer les deux concepts [Salom\'e (1998) p. 67 et p. 89]. L'inf\'erence fiduciaire sert rarement de point de vue pour juger les {\it p}-values, peut-\^etre parce qu'il n'existe pas une d\'efinition suffisamment g\'en\'erale [Buehler (1980)].

Notre regard sur les {\it p}-values provient d'une mani\`ere diff\'erente de poser le probl\`eme du choix entre deux hypoth\`eses. Nous ne cherchons pas \`a s\'election\-ner une r\`egle de d\'ecision, nous nous contentons d'un crit\`ere d\'efi\-nissant les ``bonnes'' r\`egles, celles que nous appellerons experts. Les d\'ecisions prises par ces experts sont ensuite synth\'etis\'ees dans un vote qui donne un indice de confiance pour chacune des deux hypoth\`eses. Dans un mod\`ele \`a rapport de vraisemblance monotone, la {\it p}-value d'un test unilat\'eral est un vote particulier. Pour les tests bilat\'eraux ceci n'est g\'en\'eralement possible qu'en changeant le cadre d\'ecisionnel, plus pr\'ecis\'ement la tribu du mod\`ele statistique. Avant de pr\'eciser ces r\'esultats il nous faut d\'efinir les notions d'expert et de vote. 							
\bigskip
\ESsection{Expertises}%
\medskip

L'\'etude d'un probl\`eme de d\'ecision statistique passe par la donn\'ee de crit\`eres de s\'election entre les diff\'erentes r\`egles de d\'ecision consid\'er\'ees.
Classiquement on commence par se donner une fonction de perte et on compare les proc\'edures de d\'ecision \`a partir des fonctions de risque correspondantes.
Il est rare qu'un choix unique de ce crit\`ere s'impose, bien que l'\'etude de certaines fonctions de perte soit privil\'egi\'ee, par exemple la perte quadratique dans le cadre de l'estimation ou le risque de se tromper pour le choix entre deux hypoth\`eses. M\^eme si ce dernier choix para\^{\i}t assez ``naturel", d'autres pertes sont possibles, en particulier si on regarde ce probl\`eme de choix entre deux hypoth\`eses comme un probl\`eme d'estimation [Hwang, Casella, Robert, Wells et Farrell (1992) ; Robert (1992)]. Les fonctions de risque per\-met\-tent d'introduire un pr\'eordre partiel sur les r\`egles de d\'ecision et ainsi de s\'electionner les r\`egles admissibles. Cette premi\`ere s\'election s'obtient en comparant les r\`egles entre elles, elle ne d\'erive pas d'une propri\'et\'e intrins\`eque. Nous avons cherch\'e \`a imposer ce type de propri\'et\'e pour d\'efinir les ``bonnes" r\`egles de d\'ecision, celles que nous appellerons experts. Une r\`egle pourra \^etre d\'eclar\'ee expert sans avoir \`a la comparer \`a l'ensemble des autres r\`egles. La propri\'et\'e de base que nous imposons aux experts est simplement que le dia\-gnos\-tic $d$, $\theta$ appartient \`a $\Theta_d$, doit \^etre plus probable quand $\theta$ appartient \`a $\Theta_d$ que lorsque $\theta$ n'y appartient pas. C'est une notion de r\`egle non biais\'ee, mais pour que cela suppose une connaissance fine du mod\`ele nous imposons que cette propri\'et\'e soit aussi v\'erifi\'ee conditionnellement \`a tout \'ev\'enement non n\'egligeable, et pas simplement globalement comme c'est g\'en\'eralement le cas. 
Dans le cadre de cet article, nous ne donnerons la d\'efinition des experts que pour les mod\`eles statistiques $(\Omega ,{\cal A},(\P_\theta=p_\theta.\mu)_{\theta\in\Theta})$ ayant des densit\'es strictement positives. Ceci \'evite les complications li\'ees au glissement des supports des probabilit\'es $\P_\theta$ (le cas g\'en\'eral est trait\'e dans Morel(1997)).
\medskip
\ESdefinition{D\'efinition}

{\it Consid\'erons le probl\`eme de d\'ecision d\'efini par la partition : $\Theta=\Theta_0+\Theta_1$.
La r\`egle de d\'ecision $\phi : (\Omega ,{\cal A})\rightarrow\{0,1\}$ est un expert du choix entre $\Theta_0$ et $\Theta_1$ si pour tout couple $(\theta_0,\theta_1)\in\Theta_0\times\Theta_1$ et tout \'ev\'enement non n\'egligeable $C\in{\cal A}$ elle v\'erifie :
\dli $\P_{\theta_1}(C\cap\{\phi =1\})/\P_{\theta_1}(C)\geq \P_{\theta_0}(C\cap\{\phi =1\})/\P_{\theta_0}(C)$.
}
\ESenddefinition 
\medskip 

Cette propri\'et\'e est intrins\`eque, elle ne fait pas intervenir d'autres r\`egles que $\phi$.
Dans le contexte d'une r\'eflexion sur les {\it p}-values on pourrait juger pertinent de consid\'erer les r\`egles de d\'ecision \`a valeurs dans $[0,1]$. On rajouterait ainsi de nombreuses r\`egles, par exemple toutes les fonctions croissantes du rapport $p_{\theta_1}/p_{\theta_0}$ pour le choix entre deux probabilit\'es. Il deviendrait alors difficile de d\'efinir une probabilit\'e sur l'ensemble des experts
(voir le paragraphe 2-1). Nous pouvons justifier notre d\'emarche en disant qu'au d\'epart, on cherche \`a d\'efinir des experts qui prennent totalement leurs responsabilit\'es en choisissant l'une des deux hypoth\`eses. La recherche d'un indice de confi\-ance pour les hypoth\`eses se fait dans une deuxi\`eme \'etape, si diff\'erents experts sont possibles et si l'on ne tient pas \`a essayer d'en choisir un.

Un expert du choix entre $\Theta_0$ et $\Theta_1$ est aussi un expert de tout probl\`eme de d\'ecision embo\^{\i}t\'e d\'efini par : $\Theta' = \Theta'_0 + \Theta'_1$, $\Theta'_0\subseteq\Theta_0$ et $\Theta'_1\subseteq\Theta_1$.
Le passage d'un mod\`ele embo\^{\i}t\'e au mod\`ele de d\'epart ne peut que r\'eduire l'ensemble des experts. Ceux ci peuvent donc se restreindre \`a $\II_\emptyset$ et $\II_\Omega$. Nous verrons que c'est le cas des probl\`emes bilat\`eres. Le travail sur les experts est ainsi diff\'erent de celui fait sur les r\`egles admissibles dans la th\'eorie de la d\'ecision \`a partir d'une fonction de perte. En effet, les r\`egles admissibles pour un mod\`ele embo\^{\i}t\'e sont g\'en\'eralement encore admissibles dans le mod\`ele de d\'epart. Les r\`egles admissibles sont souvent trop nombreuses, les experts eux sont plut\^ot trop rares. Ceci nous am\`ene \`a distinguer deux cas.
\medskip
\ESsubsection{Le probl\`eme de d\'ecision est expertisable}%
\medskip
Dans ce cas l'ensemble des experts ne se r\'eduit pas aux experts tri\-viaux. Pour cela il faut que les probl\`emes de d\'ecision embo\^{\i}t\'es soient semblables. Par exemple dans la th\'eorie des tests, les propri\'et\'es obtenues avec des hypoth\`eses simples (lemme de Neyman-Pearson principalement) se prolongent
facilement au cas des tests unilat\'eraux dans un mod\`ele \`a rapport de
vraisemblance monotone. Ce qui est important dans ce type de probl\`eme de d\'ecision, c'est que le probl\`eme du choix entre un \'el\'ement $\theta_0$ de
$\Theta_0$ et un \'el\'ement $\theta_1$ de $\Theta_1$ ne change pas
fondamentalement lorsque $\theta_0$ et $\theta_1$ varient. Il existe une statistique r\'eelle $T$ qui ordonne les observations, des plus favorables \`a $\theta_1$ aux plus favorables \`a $\theta_0$. Les bonnes r\`egles de d\'ecision sont de la forme $\II_{\{T<t\}}$ et $\II_{\{T\leq t\}}$, nous verrons qu'il en est de m\^eme des experts. Face \`a ce trop plein de bonnes r\`egles d\'eterministes le statisticien se donne g\'en\'eralement des crit\`eres suppl\'ementaires pour essayer de s\'electionner une r\`egle de d\'ecision : recherche d'un test U.P.P., d'un test sans biais U.P.P., d'un test invariant U.P.P., d'une r\`egle de Bayes, d'une r\`egle minimax, etc. Ceci est n\'ecessaire lorsque l'utilisateur cherche \`a mettre au point une r\`egle de d\'ecision automatique. Dans d'autres contextes c'est plut\^ot une aide \`a la d\'ecision qui est demand\'ee, un indice de confiance pour chacune des hypoth\`eses. Cette perspective nous conduit  \`a essayer de r\'esumer les choix de l'ensemble des experts. Nous allons consid\'erer nos experts comme \'egaux en droit, les faire voter et fournir \`a l'utilisateur le r\'esultat de ce vote. 
\dli Pour faire voter nos experts il faut d\'efinir une probabilit\'e sur l'ensemble des experts. Nous en avons d\'efini une \`a partir de chaque probabilit\'e $\P_\theta$ en essayant d'accorder d'autant plus de poids \`a un ensemble d'experts qu'il est form\'e d'experts donnant des r\'esultats diff\'erents sous $\P_\theta$. Par exemple, les experts $\{\II_{\{T<t\}}$, $\II_{\{T\leq t\}} ; t'<t<t''\}$ donnent la m\^eme d\'ecision sur $\{T\leq t'\}$ et $\{T\geq t''\}$, ils auront d'autant plus de poids dans le vote que $\P_\theta(\{t'<T<t''\})$ sera grand. Il est alors impossible de ``bourrer les urnes'' avec des experts presque s\^urement \'egaux. Nous avons bien s\^ur autant de votes que de valeurs du param\`etre et il y a bien des mani\`eres de les synth\'etiser. On peut faire une moyenne \`a partir d'une probabilit\'e sur l'ensemble des param\`etres en prenant en compte des informations a priori. 
On peut aussi s\'electionner le vote le plus neutre par rapport aux deux hypoth\`eses. C'est ce type de solution que nous \'etudierons ici, car il permet de donner un sens nouveau \`a la notion de {\it p}-value (la prise en compte d'informations a priori est trait\'ee dans Morel (1997)).
\medskip
\ESsubsection{Le probl\`eme de d\'ecision n'est pas expertisable}%
\medskip

Dans ce cas il n'y a que les experts tri\-viaux. 
C'est ce qui se passe g\'en\'eralement pour les hypoth\`eses bilat\'erales $H : \theta\in]-\infty,\theta_1[\cup]\theta_2,+\infty[$ et $H' : \theta\in[\theta_1,\theta_2]$. Pour obtenir un probl\`eme expertisable il faut dimi\-nuer les contraintes impos\'ees par la d\'efinition des experts en travaillant sur une sous-tribu. Ceci diminue l'ensemble des \'ev\'enements sur lesquels la propri\'et\'e de base des experts doit s'appliquer conditionnellement. Notons $H_-$ et $H_+$ les deux parties de l'hypoth\`ese englobante $H$. Si l'on ne veut pas les diff\'erencier, il ne faut pas dissocier les observations qui sont en faveur de $H_-$ de celles qui, avec la m\^eme ``force'', sont en faveur de $H_+$. Ceci revient \`a ne consid\'erer que les \'ev\'enements fournissant une information sym\'etrique sur $H_-$ et $H_+$. On retrouvera alors comme vote les {\it p}-values des tests sans biais.
\dli Cette solution ne nous semble pas la plus int\'eressante.
Dans bien des cas, le choix principal est entre $H$ et $H'$ mais $H_-$ et $H_+$ ne signifient pas la m\^eme chose. Les tests bilat\'eraux sont souvent utilis\'es comme s'ils \'etaient construits pour choisir entre les trois hypoth\`eses $H_-$,  $H'$ et $H_+$. Dans ce cas il est alors int\'eressant que le vote par rapport \`a $H$ et $H'$ soit coh\'erent avec les deux votes correspondant aux hypoth\`eses unilat\'erales $\{H_- ; H'\cup H_+\}$ et $\{H_- \cup H' ; H_+\}$ [Gabriel (1969)]. Le poids des votes en faveur de $H$ doit \^etre sup\'erieur \`a celui des votes en faveur de $H_-$ (resp. $H_+$). Nous avons pour cela introduit la notion de votes compatibles sur une famille d'hypoth\`eses unilat\'erales. Cette mani\`ere de faire a l'avantage d'imposer une coh\'erence entre les solutions de probl\`emes de d\'ecision qui reposent sur une m\^eme structuration de l'espace des param\`etres. 

Nous allons maintenant faire fonctionner cette notion d'expertise sur quelques probl\`emes de d\'ecision classiques. 

\bigskip
\ESsection{Hypoth\`eses unilat\'erales\cr\lf 
dans un mod\`ele \`a rapport de vraisemblance monotone}%
\medskip
Dans la th\'eorie des tests la premi\`ere application du lemme de Neyman-Pearson concerne les hypoth\`eses unilat\'erales dans un mod\`ele \`a rapport de vraisemblance monotone. Le cas le plus couramment trait\'e est celui d'une statistique r\'eelle exhaustive $X$ de loi $\P_\theta=p_\theta.\mu$, le param\`etre $\theta$ appartenant \`a un intervalle $\Theta$  de $\R$, la propri\'et\'e de vraisemblance monotone exprimant la croissance stricte du rapport $p_{\theta''}/p_{\theta'}$ lorsque $\theta'<\theta''$ [Monfort (1982)]. Nous nous restreindrons ici \`a des densit\'es par rapport \`a la mesure de Lebesgue $\lambda$, continues en $\theta$. Dans les autres mod\`eles les conclusions restent identiques mais nous \'eviterons les probl\`emes techniques induits par la pr\'esence de masses ou par la discontinuit\'e de la famille $\{p_\theta\}_{\theta\in\Theta}$ (le cas g\'en\'eral est trait\'e dans Morel (1997 ou 1998)). Ce cas particulier recouvre un grand nombre de cas classiques : les mod\`eles exponentiels \`a param\`etre r\'eel mais aussi les probl\`emes statistiques portant sur le param\`etre de non centralit\'e d'une famille de densit\'es de Student, Fisher ou khi-deux [Karlin (1955)].

Soit $\theta_1\in\Theta$, nous allons \'etudier le probl\`eme du choix entre $\Theta_1=]-\infty,\theta_1]$ et $\Theta_0=]\theta_1,+\infty[$ (gr\^ace \`a la continuit\'e des $p_\theta$ en $\theta_1$, prendre $\Theta_1=]-\infty,\theta_1[$ et $\Theta_0=[\theta_1,+\infty[$ ne change rien). Notons $F(\theta,x)$ la valeur en $x$ de la fonction de r\'epartition  de $X$ pour la probabilit\'e $\P_\theta$.

\dli Dans la th\'eorie des tests il y a deux traitements distincts suivant que l'on choisit de privil\'egier $\Theta_0$ ou $\Theta_1$. Le test de $H_0 : \theta\in\Theta_0$ contre $H_1 : \theta\in\Theta_1$ au seuil $\alpha$, rejette $H_0$ pour $x<t_\alpha$ avec $F(\theta_1,t_\alpha)=\alpha$. Le test de $H'_0 : \theta\in\Theta_1$ contre $H'_1 : \theta\in\Theta_0$ au seuil $\alpha$, rejette $H'_0$ pour $x>t'_\alpha$ avec $F(\theta_1,t'_\alpha)=1-\alpha$. Lorsqu'on r\'ealise $x$, les logiciels statistiques donnent pour le premier test la {\it p}-value $\alpha_{\theta_1}(x)=F(\theta_1,x)$ et pour le deuxi\`eme test la {\it p}-value $\alpha'_{\theta_1}(x)=1-F(\theta_1,x)$. Ces seuils minimum de rejet $\alpha_{\theta_1}(x)$ et $\alpha'_{\theta_1}(x)=1-\alpha_{\theta_1}(x)$ sont interpr\'et\'es comme des indices statistiques de la vraisemblance de chacune des hypoth\`eses $\Theta_0$ et $\Theta_1$. Les {\it p}-values des tests unilat\'eraux sont bien coh\'erentes, pour $\theta_1<\theta_2$ on a  $\Theta_0=]\theta_1,+\infty[\supset]\theta_2,+\infty[$ et $\alpha_{\theta_1}(x)>\alpha_{\theta_2}(x)$.
\dli Nous allons retrouver ces {\it p}-values comme votes d'experts. Les experts du choix entre $\Theta_1$ et $\Theta_0$ sont les diff\'erents tests possibles entre ces deux hypoth\`eses.
\medskip
\EStheorem{Proposition}

{\it Sur un intervalle $I\subseteq\R$ muni de la mesure de Lebesgue, consid\'erons une famille $\{p_\theta\}_{\theta\in\Theta}$ de densit\'es strictement positives. $\Theta$ \'etant un intervalle de $\R$, cette famille est suppos\'ee \`a rapport de vraisemblance strictement monotone.
Les experts du choix entre les deux hypoth\`eses unilat\'erales $\Theta_1=]-\infty,\theta_1]\cap\Theta\not=\emptyset$ et $\Theta_0=]\theta_1,+\infty[\cap\Theta\not=\emptyset$ sont les r\`egles de d\'ecision \`a valeurs dans $\{0,1\}$, presque s\^urement de la forme $f_t(x)=\II_{]-\infty,t[}(x)$ avec $t\in\overline{\R}$. 
}
\ESendtheorem 
\medskip

\ESproof{D\'emonstration}

{\it I -- $f_t=\II_{]-\infty,t[}$ est un expert.}

Les \'ev\'enements non n\'egligeables sont les bor\'eliens dont l'intersection avec l'intervalle $I$ est de mesure de Lebesgue non nulle. Soient $C$ un de ces \'ev\'enements, $\theta'\leq\theta_1$ et $\theta''>\theta_1$. On doit d\'emontrer l'in\'egalit\'e : 
\dli $\P_{\theta'}(C\cap]-\infty,t[)/\P_{\theta'}(C)\geq \P_{\theta''}(C\cap]-\infty,t[)/\P_{\theta''}(C)$. 

On va travailler avec la famille des probabilit\'es conditionnelles $\P^C_\theta$ de densit\'e $(1/\P_\theta(C))\II_C$ par rapport \`a $\P_\theta$. $f_t$ est un test de Neyman et Pearson pour tester $\P_{\theta''}^C$ contre
$\P_{\theta'}^C$ au seuil $\P_{\theta''}^C(]-\infty,t[)$. La puissance \'etant sup\'erieure au seuil [Lehmann (1986) p. 76], on a $\P_{\theta'}^C(]-\infty,t[)\geq \P_{\theta''}^C(]-\infty,t[)$, ce qui est bien l'in\'egalit\'e recherch\'ee.

\medskip
{\it II -- Un expert $\phi$ est presque s\^urement de la forme $f_t$.}

Posons $t'\,=\,sup\bigl\{t\in\overline{\R}\, ;\,f_{t}\relmont{\leq}{p.s.}\phi\bigr\}$ et  $t''\,=\,inf\bigl\{t\in\overline{\R}\, ;\,f_t\relmont{\geq}{p.s.}\phi\bigr\}$. On doit d\'emontrer que le cas $t'<t''$ est impossible lorsque $\phi$ est un expert.
Pour cela nous allons supposer $t'<t''$ et 
trouver un \'ev\'enement $C$ non n\'egligeable v\'erifiant :
$\P_{\theta'}(C\cap\{\phi =1\})/\P_{\theta'}(C)<\P_{\theta''}(C\cap\{\phi =1\})/\P_{\theta''}(C)$ lorsque $\theta'\leq\theta_1<\theta''$. 
\dli Soit $t\in]t',t''[$.
\dli Consid\'erons l'\'ev\'enement
$A=[t',t[\cap\{\phi=0\}$, par d\'efinition de $t'$ il n'est pas Lebesgue n\'egligeable. De m\^eme, la d\'efinition de $t''$ entra\^{\i}ne que l'\'ev\'enement $B=]t,t'']\cap\{\phi=1\}$ n'est pas n\'egligeable. 
\dli Posons $C=A\cup B$, on a $\P_{\theta}(C\cap\{\phi =1\})/\P_{\theta}(C)=\P_{\theta}(B)/[\P_{\theta}(A)+\P_{\theta}(B)]$. Il est facile de voir que l'\'ev\'enement $C$ v\'erifie l'in\'egalit\'e recherch\'ee si et seulement si $\P_{\theta'}(A)/\P_{\theta'}(B)>\P_{\theta''}(A)/\P_{\theta''}(B)$.
Cette derni\`ere in\'egalit\'e est une cons\'equence directe de la propri\'et\'e de vraisemblance strictement monotone ; notons $k$ le rapport $p_{\theta''}(t)/p_{\theta'}(t)$, pour $x\in A$ on a : $p_{\theta''}(x)<k p_{\theta'}(x)$ et pour $x\in B$ : $p_{\theta''}(x)>k p_{\theta'}(x)$, donc $\P_{\theta''}(A)<k \P_{\theta'}(A)$ et $\P_{\theta''}(B)>k \P_{\theta'}(B)$.
\ESendproof
\medskip

Chacun de ces experts traduit un a priori plus ou moins fort en faveur de l'une des deux hypoth\`eses. Plus $t$ est grand et plus l'expert $f_t$ parie sur l'hypoth\`ese $\theta\in\Theta_1$. Si l'on voulait en choisir un il faudrait imposer des contrain\-tes suppl\'ementaires. Dans la th\'eorie des tests, le choix de l'hypoth\`ese qui jouera le r\^ole de $H_0$ revient \`a privil\'egier les experts qui parient en faveur de $H_0$. Nous avons pr\'ef\'er\'e garder l'ensemble des avis des experts et essayer de les r\'esumer. Se pose alors le probl\`eme du poids \`a accorder \`a chacun de ces avis, il faut probabiliser l'ensemble des experts : $\{f_t\}_{t\in\overline{\R}}$. Nous quittons ici une approche de type Neyman-Wald. La difficult\'e du choix des crit\`eres de s\'election est remplac\'ee par celle de la d\'efinition du vote des experts. 
\dli Pour probabiliser l'espace des experts nous le munissons de la tribu d\'efinie \`a partir de l'ordre primordial introduit par l'indice $t\in\overline{\R}$. Nous donnerons \`a un intervalle d'experts $[f_{t'},f_{t''}]=\{f_t\, ;\, t'\leq t\leq t''\}$ d'autant plus de poids qu'il est form\'e d'experts tr\`es diff\'erents. Si $\theta$ est la vraie valeur du param\`etre il est assez naturel de juger les diff\'erences entre experts en utilisant la probabilit\'e $\P_{\theta}$. Les experts $[f_{t'},f_{t''}]$ sont \'egaux en dehors de l'intervalle $[t',t''[$ nous leur accorderons d'autant plus d'importance que la probabilit\'e $\P_\theta([t',t''[)$ est grande. Pour chaque $\theta\in\R$ nous pouvons ainsi d\'efinir sur l'ensemble des experts une probabilit\'e inductive dont la fonction de r\'epartition est $F(\theta,.)$. Lorsqu'on  r\'ealise $x$, les experts qui d\'ecident $\theta\in\Theta_1$ sont ceux correspondant \`a $t>x$. Cet ensemble d'experts $\{f_t\, ;\, t>x\}$ a alors un poids \'egal \`a $1- F(\theta,x)$. Ceci nous conduit \`a la notion de vote.

\medskip
\ESdefinition{D\'efinition}

{\it Consid\'erons le probl\`eme du choix entre les deux hypoth\`eses unilat\'erales $\Theta_1=]-\infty,\theta_1] \cap\Theta\not=\emptyset$ et $\Theta_0=]\theta_1,+\infty[\cap\Theta\not=\emptyset$, dans le mod\`ele de la proposition 3-1. On appelle vote des experts sous $\P_\theta$ la famille de probabilit\'es inductives $\{Q^x_\theta\}_{x\in\R}$ d\'efinie sur l'espace $\{0,1\}$, correspondant aux deux hypoth\`eses, par $Q^x_{\theta}(\{1\})=1-F(\theta,x)$.
}
\ESenddefinition 

\medskip

Nous avons ainsi autant de votes qu'il y a de valeurs du param\`etre. Il nous faut maintenant choisir un de ces votes ou les r\'esumer en d\'efinissant un m\'elange. Pour cela nous pouvons faire intervenir une information a priori. Les exemples d'application o\`u une telle information existe sont cependant rares. L'analyse bay\'esienne, qui repose sur la donn\'ee d'une probabilit\'e a priori, est appliqu\'ee la plupart du temps sans utiliser d'information a priori, en partant par exemple d'une loi a priori non informative. C'est d'ailleurs ce type de lois a priori qui permet de retrouver certaines {\it p}-values comme probabilit\'es a posteriori. Nous allons obtenir ces {\it p}-values comme vote d'experts en cherchant un vote ne reposant sur aucune information a priori. Nous n'\'etudierons donc pas ici la possibilit\'e de faire intervenir ce type d'information (ce cas est trait\'e dans Morel (1998)).
\dli Le vote des experts sous $\P_\theta$ est d'autant plus favorable \`a la d\'ecision $d=1$ que $\theta$ est grand. Sous $\Theta_0$ on ne veut pas trop favoriser la d\'ecision $d=1$ et inversement sous $\Theta_1$ on ne veut pas trop favoriser la d\'ecision $d=0$. Pour chacune de ces hypoth\`eses on consid\`ere le vote qui lui est le plus favorable. Sous $\Theta_1$ c'est le vote d\'efini par le point fronti\`ere $\theta_1$. Il en est de m\^eme sous $\Theta_0$ car la famille $\{p_\theta\}$ est continue en $\theta$. Le vote correspondant \`a $\theta=\theta_1$ a la propri\'et\'e d'\^etre neutre, de n'avantager aucune des deux hypoth\`eses au sens suivant : $\E_{\theta_1}[Q^x_{\theta_1}(\{1\})]=\E_{\theta_1}[Q^x_{\theta_1}(\{0\})]=1/2$ ; c'est une cons\'equence directe du fait que les valeurs de la fonction de r\'epartition se r\'epartissent uniform\'ement sur $[0,1]$. 
\dli Pour une r\'ealisation $x$, le vote au point fronti\`ere (sous $\P_{\theta_1}$) d\'efinit une probabilit\'e inductive $Q^x_{\theta_1}$ qui donne les {\it p}-values : $Q^x_{\theta_1}(\{1\})=1-F(\theta_1,x)=\alpha'_{\theta_1}(x)$ et $Q^x_{\theta_1}(\{0\})=\alpha_{\theta_1}(x)$. Ceci fournit une l\'egitimit\'e suppl\'ementaire \`a l'interpr\'etation de ces {\it p}-values comme indice de vraisemblance des hypoth\`eses unilat\'erales $\theta\leq\theta_1$ et $\theta >\theta_1$.
\medskip
\ESexample{Exemple}

En analyse de variance \`a effets fixes on utilise des statistiques qui sui\-vent des lois de Fisher d\'ecentr\'ees [Scheff\'e(1970) p. 38]. Le param\`etre de non centralit\'e $\lambda$ \'etant la quantit\'e, qui int\'eresse l'utilisateur, exprim\'ee par rapport \`a l'\'ecart-type commun des variables. Par exemple, dans une analyse de variance classique \`a deux facteurs, $\lambda$ est \'egal \`a $\sqrt{\sum(\theta_{ij})^2}/\sigma$ pour le test sur l'additivit\'e des facteurs, et \'egal \`a $\sqrt{\sum(\alpha_{i})^2}/\sigma$ pour le test sur la nullit\'e des effets additifs du premier facteur. Si une statistique $W$ suit une loi de Fisher d\'ecentr\'ee de param\`etre $\lambda$ et de degr\'es de libert\'e $(k,l)$, nous savons que $(k/l)W$ suit une loi b\^eta sur $\R^+$ : 
$\beta({k\over 2},{l\over 2},{\lambda^2\over 2})$ [Barra (1971) p. 84]. Nous allons \'etudier ce mod\`ele statistique param\'etr\'e par $\theta={\lambda^2\over 2}$. 

Consid\'erons sur $\R^+$ la famille des lois b\^eta d\'ecentr\'ees :  $\beta(p,q,\theta)$, les param\`etres $p>0$ et $q>0$ sont connus, le param\`etre de non centralit\'e $\theta\geq 0$, lui, est inconnu. Comme pour la famille des lois de Fisher d\'ecentr\'ees, on montre qu'elle est \`a rapport de vraisemblance strictement monotone [Karlin (1955)]. Soient $\theta_1\geq 0$ et une r\'ealisation $x\in\R^+$, le vote neutre en faveur de $\Theta_1=[0,\theta_1]$ est \'egal \`a : 
\cleartabs
\+ $Q^x_{\theta_1}(\{1\})$&=&$1-F(\theta_1,x)=\int_x^{+\infty}[\sum_{m=0}^{+\infty}e^{-\theta_1}.{\theta_1^m\over m!}{\Gamma(p+m+q)\over\Gamma(p+m)\Gamma(q)}{y^{p+m-1}\over (1+y)^{p+m+q}}]\,dy$\cr
\+ &=&$\int_{\N}[1- F(p+m,q,x)]\,d{\cal P}_{\theta_1}(m)$,\cr 
\dli ${\cal P}_{\theta_1}$ \'etant une loi de Poisson de param\`etre $\theta_1\geq 0$ et 
$F(p+m,q,x)$ la valeur en $x$ de la fonction de r\'epartition d'une loi $\beta(p+m,q)$. Pour $\theta_1=0$ on obtient $1-F(p,q,x)$, c'est le seuil minimum de rejet du test de $H_0 : \theta=0$ contre $H_1 : \theta>0$ (test classique en analyse de variance).
\dli La notion d'expertise apporte une justification aux {\it p}-values de l'analyse de variance. Elle permet aussi de traiter des hypoth\`eses non nulles de la forme $\theta\in[0,\theta_1]$, le r\'esultat \'etant celui donn\'e par l'inf\'erence fiduciaire.

\ESendexample

\bigskip
\ESsection{Hypoth\`eses bilat\'erales}%
\medskip
Nous allons \'etudier le probl\`eme du choix entre $\Theta_1=]-\infty,\theta_1[\cup ]\theta_2,+\infty[$ et $\Theta_0=[\theta_1,\theta_2]$ avec $\theta_1<\theta_2$ lorsque le mod\`ele statistique v\'erifie les propri\'et\'es impos\'ees dans le paragraphe pr\'ec\'edent (prendre $\Theta_0=]\theta_1,\theta_2[$ ne change rien, car les densit\'es $p_\theta$ sont suppos\'ees continues en $\theta$).
\dli Dans la th\'eorie des tests il y a encore deux traitements distincts suivant que l'on choisit de privil\'egier $\Theta_0$ ou $\Theta_1$. Pour tester $H_0 : \theta\in\Theta_0$ contre $H_1 : \theta\in\Theta_1$ il n'existe pas de test uniform\'ement plus puissant, on se restreint aux tests sans biais. 
L'\'etude de ces tests est classiquement faite dans les mod\`eles exponentiels sur $\R$ de la forme $dP_\theta(\omega)=c(\theta)exp(\theta X(\omega)).d\mu(\omega)$ [Lehmann (1986), Monfort (1982)]. Dans le cadre de cet article, seul le cas o\`u la mesure $\mu$ est diffuse nous int\'eresse.
Au seuil $\alpha$ le plus puissant des tests sans biais rejette alors $H_0$ pour $x<c_1$ et $x>c_2$, $c_1$ et $c_2$ \'etant d\'etermin\'es par $F(\theta_1,c_1)+[1-F(\theta_1,c_2)]= F(\theta_2,c_1)+[1-F(\theta_2,c_2)]=\alpha$. Pour tester $H'_0 : \theta\in\Theta_1$ contre $H'_1 : \theta\in\Theta_0$ au seuil $\alpha$ il existe un test uniform\'ement plus puissant. Il rejette $H'_0$ lorsque $c'_1\leq x\leq c'_2$, $c'_1$ et $c'_2$ \'etant d\'etermin\'es par $F(\theta_1,c_2)-F(\theta_1,c_1)= F(\theta_2,c_2)-F(\theta_2,c_1)=\alpha$ . 
\dli Les seuils minimum de rejet fournis par ces tests sont souvent interpr\'et\'es comme des indices statistiques de la vraisemblance de chacune des hypoth\`eses $\Theta_0$ et $\Theta_1$. Nous avons vu dans l'introduction que cette mani\`ere de faire est soumise \`a de nombreuses critiques. Nous allons montrer que dans le cadre d'une expertise elle n'est justifi\'ee qu'en changeant le mod\`ele statistique. 
\dli Tel qu'il est pos\'e, le probl\`eme de d\'ecision entre deux hypoth\`eses bilat\'erales ne poss\`ede que $\II_\emptyset$ et $\II_{\R}$ comme experts. Ceci est du au fait que tout expert du choix entre deux hypoth\`eses est aussi un expert des probl\`emes de d\'ecision embo\^{\i}t\'es, en particulier ici du choix entre les hypoth\`eses unilat\'erales $]-\infty,\theta_1[$ et $\Theta_0=[\theta_1,\theta_2]$, mais aussi du choix entre $\Theta_0$ et $]\theta_2,+\infty[$. C'est le m\^eme argument qui montre qu'il n'y a pas de test uniform\'ement plus puissant de $H_0 : \theta\in\Theta_0$ contre $H_1 : \theta\in\Theta_1$. 
\dli Face \`a cette p\'enurie d'experts nous envisagerons deux traitements possibles. On peut vouloir une coh\'erence entre la solution \`a ce probl\`eme et les expertises des hypoth\`eses unilat\'erales d\'eriv\'ees. Cette mani\`ere de faire permet de conserver le sens donn\'e par l'ordre qui structure l'espace des param\`etres. On peut aussi envisager de casser cet ordre en cr\'eant une  \'equivalence entre les \'el\'ements des deux parties de $\Theta_1$ : $\Theta^-_1=]-\infty,\theta_1[$ et $\Theta^+_1=]\theta_2,+\infty[$, c'est ce que fait la contrainte de non-biais impos\'ee dans la th\'eorie des tests de $\Theta_0$ contre $\Theta_1$.
\medskip
\ESsubsection{Expertises compatibles}%
\medskip
Bien souvent le choix entre $\Theta_1$ et $\Theta_0$ cache un probl\`eme de d\'ecision entre trois hypoth\`eses : $\Theta^-_1$, $\Theta_0$ et $\Theta^+_1$.
L'utilisateur veut d'abord savoir si $\theta$ appartient \`a $\Theta_1$ ou $\Theta_0$, mais le regroupement de $\Theta^-_1$ et $\Theta^+_1$ peut cacher des interpr\'etations de type diff\'erent suivant que $\theta$ est inf\'erieur \`a $\theta_1$ ou sup\'erieur \`a $\theta_2$. Par exemple, si $\Theta_0=[-\theta_0,+\theta_0]$ est la traduction d'un effet n\'egligeable, le rejet de cette hypoth\`ese conduit g\'en\'eralement dans la pratique \`a des conclusions diff\'erentes suivant que l'estimation de $\theta$ est positive ou n\'egative. Dans un tel cadre le choix entre les hypoth\`eses unilat\'erales $\Theta^-_1$ et $\Theta_0\cup\Theta^+_1$ est envisageable, ainsi que le choix entre $\Theta^-_1\cup\Theta_0$ et $\Theta^+_1$. Il est alors primordial que les trois probl\`emes de d\'ecision g\'en\'er\'es par $\Theta^-_1$, $\Theta_0$ et $\Theta^+_1$ aient des solutions coh\'erentes. Les votes s\'electionn\'es dans les deux probl\`emes unilat\'eraux doivent alors conduire \`a un vote en faveur de $\Theta^-_1$ plus faible que celui en faveur de $\Theta^-_1\cup\Theta_0$. Nous dirons qu'ils sont compa\-ti\-bles, ils permettent alors d'induire une probabilit\'e sur la tribu engendr\'ee par les trois hypoth\`eses. Les trois probl\`emes de d\'ecision consid\'er\'es ne sont g\'en\'eralement pas les seuls envisageables dans une application donn\'ee. Il est en effet rare que les bornes $\theta_1$ et $\theta_2$ ne soient pas sujettes \`a discussion. Aussi nous avons \'etendu la notion de votes compatibles \`a un ensemble quelconque d'hypoth\`eses unilat\'erales.

\medskip
\ESdefinition{D\'efinition}

{\it Soient $(\{\Theta_1^f,\Theta_0^f\})_{f\in{\cal F}}$ une famille d'hypoth\`eses unilat\'e\-rales et $Q$ une application de $\R\times\{\Theta_1^f\}_{f\in{\cal F}}$ dans $[0,1]$ telle que $Q^x(\Theta_1^{f})$ repr\'e\-sente la valeur en $d=1$ d'un vote d'experts du choix entre $\Theta_1^{f}$ et $\Theta_0^{f}$ lorsqu'on r\'ealise $x$.
\dli Nous dirons que $Q$ d\'efinit des votes compatibles lorsqu'elle se prolonge, pour presque tout $x$, en une probabilit\'e unique $Q^x$ sur la tribu engendr\'ee par $\{\Theta_1^f\}_{f\in{\cal F}}$.
}
\ESenddefinition 
\medskip
Dans le cadre de cet article nous n'\'etudierons la compatibilit\'e des votes que pour les votes neutres s\'electionn\'es au paragraphe pr\'ec\'edent (l'intro\-duction d'informations a priori est trait\'ee dans Morel (1997)).
Pour les hypoth\`eses unilat\'erales $\{]-\infty,\theta],]\theta,+\infty[\}$ le vote neutre est celui d\'efini par le point fronti\`ere $\theta$. Nous allons donc consid\'erer l'application : 

\dli $Q^x(]-\infty,\theta])=1-F(\theta,x)$, $F(\theta,.)$ \'etant la fonction de r\'epartition de $\P_\theta$. 
\medskip
\EStheorem{Proposition}
{\it Sur un intervalle $I\subseteq\R$ muni de la mesure de Lebesgue, consid\'erons une famille $\{p_\theta\}_{\theta\in\Theta}$ de densit\'es strictement positives. $\Theta$ \'etant un intervalle de $\R$, cette famille est suppos\'ee continue en $\theta$ et \`a rapport de vraisemblance strictement monotone.
\dli Pour tout $\theta\in\Theta-\{sup\Theta\}$ consid\'erons le vote neutre d\'efini par les hypoth\`eses unilat\'erales $\Theta_1^\theta =]-\infty,\theta]\cap\Theta $ et $\Theta_0^\theta =]\theta,+\infty[\cap\Theta$. 
Ces votes sont compatibles si et seulement si, pour presque tout $x\in I$ on a : 
\dli 1) $lim_{\theta\rightarrow sup\Theta}F(\theta,x)=0$ si $sup\Theta\notin\Theta$
\dli 2) $lim_{\theta\rightarrow inf\Theta}F(\theta,x)=1$ si $inf\Theta\notin\Theta$. 
}
\ESendtheorem 
\medskip

\ESproof{D\'emonstration}

Les votes neutres d\'efinissent une application $Q$ de $\R\times\{\Theta_1^\theta\}_{\theta\in\Theta-{sup\Theta}}$ dans $[0,1]$ : $Q^x(\Theta_1^{\theta})=1-F(\theta ,x)$. 
$Q$ donne des votes compatibles si elle se prolonge, pour presque tout $x$, en une probabilit\'e unique sur les Bor\'eliens de l'intervalle $\Theta$. 
\dli Pour une r\'ealisation $x$ donn\'ee, l'application $Q^x$ est prolongeable si elle v\'erifie les propri\'et\'es caract\'eristiques d'une fonction de r\'epartition.
\dli a) La croissance en $\theta$ de $Q^x(\Theta_1^{\theta})=1-F(\theta,x)$ est une cons\'equence directe de la d\'ecroissance de $F(.,x)$ sur une famille \`a rapport de vraisemblance monotone.
\dli b) La continuit\'e \`a droite en $\theta$ de $Q^x(\Theta_1^{\theta})$ est \'equivalente \`a celle de $F(\theta,x)$. La famille des densit\'es $p_\theta$ \'etant continue, elle l'est aussi dans $L_1$, ce qui implique la continuit\'e de $F(.,x)$. 
\dli c) Il reste \`a v\'erifier les conditions  aux bornes, c'est-\`a-dire la convergence de $Q^x(\Theta_1^{\theta})$ vers $1$ (resp. $0$) lorsque $\Theta_1^{\theta}$ cro\^{\i}t vers $\Theta$ (resp. d\'ecro\^{\i}t vers $\emptyset$). Une telle situation n'est possible que si $sup\Theta\notin\Theta$ (resp. $inf\Theta\notin\Theta$). Les deux conditions n\'ecessaires et suffisantes de la proposition en d\'ecoulent imm\'ediatement. 
\ESendproof
\medskip

Un ensemble de votes compatibles sur les hypoth\`eses unilat\'erales permet de prolonger de fa\c con coh\'erente ces expertises \`a des hypoth\`eses non expertisables. Sous les conditions de la proposition pr\'ec\'edente, on obtient pour les hypoth\`eses bilat\'erales $\Theta_1=]-\infty,\theta_1[\cup ]\theta_2,+\infty[$ et $\Theta_0=[\theta_1,\theta_2]$ les votes suivants : 
$Q^x(\Theta_1)=1-F(\theta_1,x)+F(\theta_2,x)$ et $Q^x(\Theta_0)=F(\theta_1,x)-F(\theta_2,x)$.
\dli Ces votes ne sont pas les {\it p}-values donn\'ees par les tests classiques. Ils nous semblent pr\'ef\'erables lorsque la structure d'ordre sur $\Theta$ reste primordiale pour l'interpr\'etation, malgr\'e la cassure introduite par l'hypoth\`ese $\Theta_1$ qui regroupe $]-\infty,\theta_1[$ et  $]\theta_2,+\infty[$. On \'evite ainsi les incoh\'erences des {\it p}-values relev\'ees par Schervish (1996). Les votes obtenus sans information a priori d\'efinissent sur $\Theta$ la
distribution fiduciaire de Fisher.
\medskip
\ESexample{Exemple}

Consid\'erons le cas trait\'e par Schervish (1996), celui de la famille des lois normales $\{N(\theta,1)\}_{\theta\in\R}$. La proposition 4.2 s'applique : $F(\theta,x)=\Phi(x-\theta)$, $\Phi$ \'etant la fonction de r\'epartition de la loi $N(0,1)$. La probabilit\'e $Q^x$ sur $\Theta$ est la loi $N(x,1)$. C'est la loi a posteriori correspondant \`a la mesure de Lebesgue, loi a priori impropre unanimement consid\'er\'ee comme non informative. 

La probabilit\'e inductive $Q^x$ obtenue \`a partir de votes neutres compatibles correspond souvent \`a la loi a posteriori d'une loi a priori qui peut \^etre consid\'er\'ee comme non informative [Morel(1997)]. Dans le cadre de ce paragraphe citons la famille des lois normales $\{N(m,\theta^2)\}_{\theta\in\R^+_*}$ et celle des lois gamma $\{\gamma(p,\theta)\}_{\theta\in\R_*^+}$ (voir l'exemple 5.1).
\ESendexample

\medskip
\ESsubsection{{\it P}-values des tests sans biais comme votes d'experts}%
\medskip
Nous consid\'erons maintenant le cas o\`u les deux parties $\Theta^-_1=]-\infty,\theta_1[$ et $\Theta^+_1=]\theta_2,+\infty[$ de $\Theta_1$ sont \'equivalentes pour l'interpr\'etation. Celle-ci est cens\'ee \'evoluer de la m\^eme mani\`ere de $-\infty$ \`a $\theta_1$ et de $+\infty$ \`a $\theta_2$. 

Dans l'exemple simple de la famille des lois normales $\{N(\theta,1)\}_{\theta\in\R}$, le probl\`eme de d\'ecision peut \^etre consid\'er\'e comme invariant par sym\'etrie autour de $c=(\theta_1+\theta_2)/2$. Cette notion sert classiquement \`a diminuer l'ensemble des r\`egles possibles. Ici nous devons au contraire augmenter le nombre d'experts. Pour cela nous allons restreindre
la tribu du mod\`ele aux \'ev\'enements invariants, les bor\'eliens sym\'etriques par rapport \`a $c$. Cette restriction sur les \'ev\'enements impose des contraintes moins fortes pour le label expert. Elle permet d'avoir d'autres experts que les experts triviaux du mod\`ele de d\'epart. En fait on travaille sur le mod\`ele associ\'e \`a la statistique $W(x)=(x-c)^2$. C'est exactement ce que l'on fait, quand pour comparer deux moyennes, on utilise la statistique $W=T^2$ de Fisher au lieu de la statistique $T$ de Student. Dans le nouveau mod\`ele $\theta$ et $2c-\theta$ d\'efinissent la m\^eme probabilit\'e, la loi de khi-deux d\'ecentr\'ee \`a un degr\'e de libert\'e et de param\`etre d'excentricit\'e $\lambda^2=(\theta-c)^2$. On peut prendre comme ensemble des param\`etres $\Lambda=\R^+$, le choix entre $\Theta_0$ et $\Theta_1$ devient un choix entre $[0,(\theta_2-\theta_1)/2]$ et $](\theta_2-\theta_1)/2,+\infty[$. Ce sont des hypoth\`eses unilat\'erales dans un mod\`ele \`a rapport de vraisemblance monotone [Karlin (1955)].
La proposition 3.1, donne les experts de ce nouveau probl\`eme de d\'ecision param\'etr\'e par $\lambda$. Avec le param\`etre $\theta$ on a le r\'esultat suivant.

\medskip
\EStheorem{Corollaire}

{\it Consid\'erons le mod\`ele statistique 
$(\R,{\cal C},\{N(\theta,1)\}_{\theta\in\R})$, ${\cal C}$ \'etant la tribu des bor\'eliens sym\'etriques par rapport \`a $c\in\R$. Posons $\theta_1=c-\lambda_1$ et $\theta_2=c+\lambda_1$ avec $\lambda_1\geq 0$. 
Les experts du choix entre les deux hypoth\`eses $\Theta_1=]-\infty,\theta_1[\cup]\theta_2,+\infty[$ et $\Theta_0=[\theta_1,\theta_2]$ sont les r\`egles de d\'ecision \`a valeurs dans $\{0,1\}$, presque s\^urement de la forme :  \dli $g_t(x)=\II_{[t,+\infty[}(\mid\! x-c\!\mid)$ avec $t\in\overline{\R^+}$. 
}
\ESendtheorem
\medskip
Les hypoth\`eses unilat\'erales $\Lambda_1=[0,\lambda_1]$ et $\Lambda_0=]\lambda_1,+\infty[$ sont \'equivalentes aux hypoth\`eses $\Theta_0$ et $\Theta_1$. Comme au paragraphe 3 on peut d\'efinir le vote au point fronti\`ere $\lambda_1$ ; lorsqu'on r\'ealise $x$, le vote en faveur de $\Lambda_0$, c'est-\`a-dire $\Theta_1$, est alors donn\'e par la valeur en $W(x)=(x-c)^2$ de la fonction de r\'epartition d'une loi de 
khi-deux d\'ecentr\'ee \`a un degr\'e de libert\'e et de param\`etre d'excentricit\'e $\lambda_1^2=(\theta_1-c)^2$. Notons $Q^x(\Theta_1)$ ce vote, c'est la probabilit\'e de l'intervalle 
$]-\mid\! x-c\!\mid,\mid\! x-c\!\mid[$  pour la loi $N(\theta_1-c,1)$. 
\dli $Q^x(\Theta_1)=\Phi(\mid\! x-c\!\mid -\theta_1+c)-\Phi(-\mid\! x-c\!\mid -\theta_1+c)$ est donc la {\it p}-value du test de $H'_0 : \theta\in\Theta_1$ contre $H'_1 : \theta\in\Theta_0$. Le vote $Q^x(\Theta_0)$ en faveur de $\Theta_0$ est bien s\^ur la {\it p}-value du test sans biais de $H_0 : \theta\in\Theta_0$ contre $H_1 : \theta\in\Theta_1$. Ceci reste valable pour le cas particulier $\theta_1=\theta_2=c$ qui correspond \`a l'hypoth\`ese $\Theta_0=[c]$. La probabilit\'e sur $\Lambda=\R^+$ donn\'ee par la proposition 4.2 poss\`ede une masse en $\{c\}$ qui est la {\it p}-value du test de $H_0 : \theta=c$ contre $H_1 : \theta\not= c $. 
\dli Sur cet exemple, nous venons de retrouver comme votes d'experts les {\it p}-values des tests bilat\'eraux classiques. Ceci n'a \'et\'e possible qu'en supposant que l'interpr\'etation du param\`etre pr\'esente une sym\'etrie par rapport au centre de l'intervalle $[\theta_1,\theta_2]$. Ce n'est pas le cas lorsque l'interpr\'etation reste li\'ee \`a l'ordre sur $\Theta$, il est alors pr\'ef\'erable d'utiliser les votes compatibles du paragraphe pr\'ec\'edent.
\dli Ce que nous venons de faire sur l'exemple simple des lois $\{N(\theta,1)\}$ peut se g\'en\'eraliser aux autres mod\`eles exponentiels. La notion de sans biais pour le test de $H_0 : \theta\in\Theta_0$ contre $H_1 : \theta\in\Theta_1$ et celle de seuil pour le test de $H'_0 : \theta\in\Theta_1$ contre $H'_1 : \theta\in\Theta_0$ imposent des r\`egles de d\'ecision dont la r\'egion de rejet $C$ v\'erifie $\P_{\theta_1}(C)=\P_{\theta_2}(C)$. Cette restriction sur les r\`egles possibles est une mani\`ere d'imposer un traitement ``sym\'etrique'' des deux parties $\Theta_1^-$ et $\Theta_1^+$ de l'hypoth\`ese $\Theta_1$. Pour retrouver les {\it p}-values de ces tests comme votes d'experts, il faut rendre le probl\`eme de d\'ecision expertisable en enlevant des contraintes dans la d\'efinition 2.1. La recherche d'une sym\'etrie de traitement des sous hypoth\`eses $\Theta_1^-$ et $\Theta_1^+$ se traduira par la restriction aux \'ev\'enements $C$ v\'erifiant $\P_{\theta_1}(C)=\P_{\theta_2}(C)$. 

\bigskip
\ESsection{Mod\`eles avec param\`etre fant\^ome}%
\medskip
L'espace des param\`etres est un espace produit $\Theta\times\Upsilon$, $\theta\in\Theta\subseteq\R$ est le param\`etre d'int\'er\^et et $\upsilon\in\Upsilon$ le param\`etre fant\^ome. On consid\`ere le probl\`eme du choix entre $\Theta_1\times\Upsilon$ et $\Theta_0\times\Upsilon$. Pour chaque valeur possible $\upsilon$ du param\`etre fant\^ome le probl\`eme de d\'ecision se r\'eduit au type de ceux que nous avons \'etudi\'es pr\'ec\'edemment. Une expertise, si elle existe, donnera des votes, en faveur de chacune des hypoth\`eses $\Theta_1$ et $\Theta_0$, index\'es par $\upsilon\in\Upsilon$. Ces r\'esultats conditionnels au param\`etre fant\^ome peuvent \^etre r\'esum\'es en uti\-li\-sant une probabilit\'e sur $\Upsilon$. Pour d\'efinir cette probabilit\'e nous chercherons \`a utiliser des votes compatibles.

Nous allons donner des exemples dans le cas simple mais courant o\`u le probl\`eme de d\'ecision repose sur deux statistiques ind\'ependantes r\'eelles $(T,U)$, qui d\'efinissent chacune un mod\`ele \`a rapport de vraisemblance monotone par rapport \`a $\theta$ respectivement $\upsilon$, la loi de $T$ pouvant d\'ependre de $\upsilon$ (d'autres cas sont trait\'es dans Morel (1997)).
\medskip
\ESexample{Exemple}

Commen\c cons par \'etudier le cas classique d'un n-\'echantillon d'une loi $N(\theta,\upsilon=\sigma^2)$, $n\geq 2$, $\theta\in\R$ et $\upsilon\in\R_*^+$.
La moyenne et la variance empiriques, $\overline{X}$ et $S^2={1\over n-1}\sum_{i=1}^n(X_i-\overline{X})^2$, sont exhaustives et ind\'ependantes. La statistique $T=\overline{X}$ est de loi $N(\theta,\upsilon/n)$, elle d\'efinit bien un mod\`ele \`a rapport de vraisemblance monotone par rapport \`a $\theta$. Quant \`a la statistique $U=(n-1)S^2$ elle suit une loi de khi-deux \`a $(n-1)$ degr\'es de libert\'e et ayant pour param\`etre d'\'echelle $\upsilon$, c'est-\`a-dire une loi gamma $\gamma({n-1\over 2},2\upsilon)$. Elle d\'efinit un mod\`ele statistique \`a rapport de vraisemblance strictement monotone par rapport \`a $\upsilon$ [Karlin (1955)].

Consid\'erons les hypoth\`eses unilat\'erales : $\Theta_1^\theta=]-\infty,\theta]$ et $\Theta_0^\theta=]\theta,+\infty[$. 
\dli Pour toute valeur fix\'ee $\upsilon$ du param\`etre fant\^ome, nous savons expertiser ce probl\`eme de d\'ecision. Il d\'epend uniquement de la statistique $T$. La proposition 4.2 s'applique et pour une r\'ealisation $t=\overline{x}$ les votes neutres compatibles d\'efinissent sur $\Theta=\R$ la loi $N(\overline{x},\upsilon/n)$.
\dli La statistique $U$ ind\'ependante de $T$ va nous permettre de probabiliser $\Upsilon=\R_*^+$. En effet elle d\'efinit sur $\R^+$ muni de la mesure de Lebesgue la famille de densit\'es $\{{1\over\Gamma(p)2\upsilon}({u\over2\upsilon})^{p-1}exp(-{u\over2\upsilon})\}_{\upsilon\in\R_*^+}$, avec $p={n-1\over 2}$. Ces densit\'es sont nulles uniquement en $u=0$ qui est n\'egligeable et la fonction de r\'epartition est donn\'ee par $F(\upsilon,u)=\int_0^{u\over2\upsilon}{1\over\Gamma(p)}y^{p-1}e^{-y}\,dy$, on peut donc appliquer la proposition 4.2. Elle nous donne les votes compatibles suivant :
\dli $Q^u(]0,\upsilon])=1-F(\upsilon,u)=\int_0^\upsilon{2\over\Gamma(p)u}({u\over2\lambda})^{p+1}exp(-{u\over2\lambda})\,d\lambda$ ($u>0$).
\dli Leur prolongement est donc la loi inverse d'une loi 
$\gamma(p,{2\over u})$, not\'ee $Inv\gamma(p,{2\over u})$. C'est la loi a posteriori pour la loi a priori impropre et non informative de densit\'e ${1\over\upsilon}$ [Berger (1985)], mais aussi la distribution fiduciaire [Fisher (1935)]. 
\dli Avec cette probabilit\'e inductive nous devons maintenant faire la moyenne des votes obtenus conditionnellement \`a $\upsilon$. Si on r\'ealise $t=\overline{x}$ et $u=(n-1)s^2$, le vote en faveur de $\Theta_1^\theta=]-\infty,\theta]$, not\'e $Q^{(t,u)}(]-\infty,\theta])$, s'obtient en faisant la moyenne de $\Phi({\theta-t\over\sqrt{\upsilon/n}})$, $\Phi$ \'etant la fonction de r\'epartition de la loi $N(0,1)$.
\medskip
\cleartabs
\+ $Q^{(t,u)}(]\!-\!\infty,\theta])$&=&$\int\Phi((\theta-t)/\sqrt{\upsilon/n})\,d\,Inv\gamma(p,{2\over u})(\upsilon)$\cr
\+ &=&$\int_{\R^+}\Phi(\sqrt{n}(\theta-t)\sqrt{\lambda})\,d\,\gamma(p,{2\over u})(\lambda)$\cr
\+ &=&$\int_{\R^+}\Phi(\sqrt{n}(\theta-t)\sqrt{\lambda}){u\over2\Gamma(p)}({u\over2}\lambda)^{p-1}\,e^{-{u\over2}\lambda}\,d\lambda$\cr
\+ &=&$\int_{\R^+}[\int_{-\infty}^{\sqrt{n}(\theta-t)}{\sqrt{\lambda}\over\sqrt{2\pi}}e^{-{x^2\over2}\lambda}\,dx]{u\over2\Gamma(p)}({u\over2}\lambda)^{p-1}\,e^{-{u \over2}\lambda}\,d\lambda$\cr
\+ &=&$\int_{\R^+}[\int_{-\infty}^{\sqrt{n}(\theta-t)\over\sqrt{s^2}}{\sqrt{s^2\lambda}\over\sqrt{2\pi}}e^{-{y^2\over2}s^2\lambda}\,dy]{u\over2\Gamma(p)}({u\over2}\lambda)^{p-1}\,e^{-{u \over2}\lambda}\,d\lambda$\cr 
\+ &=&$\int_{-\infty}^{\sqrt{n}(\theta-t)\over\sqrt{s^2}}[\int_0^{+\infty}{\sqrt{\nu}\over\sqrt{2\pi}}e^{-{y^2\over2}\nu}\times{({n-1\over 2})^{n-1\over 2}\over\Gamma({n-1\over 2})}\nu^{{n-1\over 2}-1}\,e^{-{n-1\over2}\nu}\,d\nu]\,dy$\cr
\dli L'int\'egrale entre crochets donne la densit\'e d'une loi de Student \`a $(n-1)$ degr\'es de libert\'e comme m\'elange des lois $N(0,{1\over\nu})$ par la loi $\gamma({n-1\over2},{2\over n-1})$ [Dickey (1968)].
Le vote $Q^{(t,u)}(]-\infty,\theta])$ est donc \'egal \`a la valeur de la fonction de r\'epartition d'un Student $(n-1)$ en ${\sqrt{n}(\theta-t)\over\sqrt{s^2}}=-\sqrt{n}{(\overline{x}-\theta)\over\sqrt{s^2}}$. C'est le seuil minimum de rejet du test de Student de $\Theta_1^\theta=]-\infty,\theta]$ contre $\Theta_0^\theta=]\theta,+\infty[$. Bien s\^ur
$Q^{(t,u)}(\Theta_0^\theta)=1-Q^{(t,u)}(\Theta_1^\theta)$ est le seuil minimum de rejet du test de Student de $\Theta_0^\theta$ contre $\Theta_1^\theta$.
\dli Lorsque $\theta$ parcourt $\R$, les votes pr\'ec\'edents sont compatibles. Ils d\'efinissent, sur l'espace $\R$ du param\`etre $\theta$, une probabilit\'e qui est une loi de Student \`a $(n-1)$ degr\'es de libert\'e, de moyenne $\overline{x}$ et de param\`etre d'\'echelle ${s^2\over n}$ [Dickey (1968)], c'est la distribution fiduciaire donn\'ee par Fisher (1935).

Consid\'erons maintenant le cas des hypoth\`eses bilat\'erales 
\dli $\Theta_1=]-\infty,\theta_1[\cup]\theta_2,+\infty[$ et $\Theta_0=[\theta_1,\theta_2]$. Si la structure d'ordre sur $\Theta=\R$ est primordiale pour l'interpr\'etation, il nous semble normal de prendre un vote coh\'erent par rapport aux votes correspondant aux hypoth\`eses unilat\'erales. Avec les votes neutres on obtient la
probabilit\'e inductive de $\Theta_1$ et $\Theta_0$ donn\'ee par la loi pr\'ec\'edente. Ce ne sont pas les {\it p}-values des tests de Student bilat\'eraux. Comme dans le cas des lois $N(\theta,1)$ on retrouve ces {\it p}-values en utilisant la sym\'etrie par rapport \`a $c=(\theta_1+\theta_2)/2$  du probl\`eme de d\'ecision pour le rendre expertisable. C'est en fait travailler sur $\theta$ \`a  partir des statistiques ind\'ependantes $(\overline{X}-c)^2$ et $S^2$. Nous allons le faire dans un cadre plus g\'en\'eral, celui de l'analyse de variance.

\ESendexample

\medskip
\ESexample{Exemple}

Reprenons l'exemple 3.3 de l'analyse de variance \`a effets fixes. Nous avons travaill\'e sur le param\`etre de non centralit\'e $\lambda$ d'une loi de Fisher d\'ecentr\'ee. $\lambda$ d\'epend du param\`etre fant\^ome $\sigma^2$, la variance inconnue. On a $\lambda^2=\theta/\sigma^2$, le param\`etre $\sqrt{\theta}$ est en fait g\'en\'eralement celui qui int\'eresse l'utilisateur, son expression $\lambda$ en unit\'e d'\'ecart-type est un produit de la statistique.
D'ailleurs pour les tests courants de comparaison de $\lambda$ \`a $0$, on exprime $H_0 : \lambda=0$ par $\theta=0$. Le param\`etre $\lambda$ n'appara\^{\i}t que dans la puissance de ces tests.

La statistique de Fisher d\'ecentr\'ee $W$ que nous avons utilis\'ee pr\'ec\'e\-dem\-ment est construite \`a partir de deux statistiques ind\'ependantes $T$ et $U$ ; $T/\sigma^2$ suit une loi de khi-deux d\'ecentr\'ee \`a $k$ degr\'es de libert\'e et de param\`etre de non centralit\'e $\theta$ ;  $U/\sigma^2$ suit une loi de khi-deux \`a $l$ degr\'es de libert\'e [Scheff\'e (1970)]. La r\'eduction du probl\`eme de base \`a l'\'etude du mod\`ele engendr\'e par $(T,U)$ se fait en imposant des propri\'et\'es d'invariance. Le mod\`ele image s'\'ecrit 
$((\R^+)^2,{\cal B}^2,(P^T_{(\theta,\sigma^2)}\otimes P^U_{\sigma^2})_{(\theta,\sigma^2)\in\R^+\times\R^+_*})$, $P^T_{(\theta,\sigma^2)}$ est une loi $\gamma({k\over 2},2\sigma^2)$ d\'ecentr\'ee de $\theta/(2\sigma^2)$ et $P^U_{\sigma^2}$ une loi 
$\gamma({l\over 2},2\sigma^2)$.

Plus g\'en\'eralement nous allons \'etudier les mod\`eles statistiques de la forme : 
$((\R^+)^2,{\cal B}^2,(P^T_{(\theta,\upsilon)}\otimes P^U_{\upsilon})_{(\theta,\upsilon)\in\R^+\times\R^+_*})$ o\`u $P^T_{(\theta,\upsilon)}$ est une loi $\gamma(p,\upsilon)$ d\'ecentr\'ee de $\theta/\upsilon$ et $P^U_{\upsilon}$ une loi $\gamma(q,\upsilon)$. Nous allons proc\'eder comme dans l'exemple pr\'ec\'edent. 

Consid\'erons les hypoth\`eses unilat\'erales : 
$\Theta_1^\theta=[0,\theta]$ et $\Theta_0^\theta=]\theta,+\infty[$  ($\theta\geq 0$). 
\dli Pour $\upsilon$ fix\'e le probl\`eme de d\'ecision ne d\'epend que de la statistique $T$ donc des lois  
$\{P^T_{(\theta,\upsilon)}\}_{\theta\in\R^+}$ ou de la famille des densit\'es 
\dli $\{\int_{\N}{1\over\Gamma(p+m) \upsilon}\,({t\over\upsilon})^{p+m-1}\, exp(-{t\over\upsilon})\,d{\cal P}_{\theta\over\upsilon}(m)\}_{\theta\in\R^+}$, ${\cal P}_{\theta\over\upsilon}$ \'etant une loi de poisson de param\`etre ${\theta\over\upsilon}\geq 0$ [Barra (1971)]. On est ramen\'e au choix entre deux hypoth\`eses unilat\'erales dans un mod\`ele \`a rapport de vraisemblance strictement monotone [Karlin (1955)]. Si l'on ne veut ou peut pas faire intervenir une information a priori on utilise le vote neutre, celui au point fronti\`ere $\theta$. Pour la r\'ealisation $t$ la fr\'equence des experts en faveur de $\Theta_1^\theta=[0,\theta]$ est alors \'egale \`a $1-F_\upsilon(\theta,t)$, $F_\upsilon(\theta,t)$ \'etant la valeur en $t$ de la fonction de r\'epartition de $P^T_{(\theta,\upsilon)}$.
\dli $F_\upsilon(\theta,t)=\int_{\N}\Gamma(p+m,\upsilon,t)\,d{\cal P}_{\theta\over\upsilon}(m)$, le terme $\Gamma(p+m,\upsilon,t)$ d\'esignant la valeur de la fonction de r\'epartition d'une loi $\gamma(p+m,\upsilon)$ en $t$ : \dli $\Gamma(p+m,\upsilon,t)=\int_0^t{1\over\Gamma(p+m) \upsilon}\,({x\over\upsilon})^{p+m-1}\, exp(-{x\over\upsilon})\,dx=\Gamma(p+m,1,{t\over\upsilon})$.
\dli Nous allons maintenant probabiliser $\Upsilon=\R^+_*$ en utilisant la statistique $U$ qui est de loi $\gamma(q,\upsilon)$. Nous avons \'etudi\'e ce type de mod\`ele dans l'exemple pr\'ec\'edent. Les votes neutres pour les hypoth\`eses unilat\'erales sont compa\-ti\-bles. Ils se prolongent en une probabilit\'e sur les bor\'eliens de $\Upsilon=\R^+_*$ qui est l'inverse d'une loi $\gamma(q,{1\over u})$ ($u>0$).
\dli Pour toute r\'ealisation $(t,u)$, la moyenne $Q^{(t,u)}([0,\theta])$ des votes en faveur de $[0,\theta]$ est alors \'egale \`a : 
\medskip
\cleartabs
\+ $Q^{(t,u)}([0,\theta])$&=&$1-\int_{\R^+_*}F_\upsilon(\theta,t){1\over\Gamma(q) u}\,({u\over\upsilon})^{q+1}\, exp(-{u\over\upsilon})\,d\upsilon$ \cr
\+ &=&$1-\int_{\R^+_*}\int_{\N}[\int_0^t{1\over\Gamma(p+m) \upsilon}\,({x\over\upsilon})^{p+m-1}\, exp(-{x\over\upsilon})\,dx]\,d{\cal P}_{\theta\over\upsilon}(m)\times$& \cr
\+ & &\qquad\hfill ${1\over\Gamma(q) u}\,({u\over\upsilon})^{q+1}\, exp(-{u\over\upsilon})\,d\upsilon$& \cr
\+ &=&$1-\sum_{m=0}^{+\infty}\int_0^t[\int_{\R^+_*}{exp(-{\theta\over\upsilon})\,({\theta\over\upsilon})^m\over m!\Gamma(p+m) \upsilon}\,({x\over\upsilon})^{p+m-1}\, exp(-{x\over\upsilon})\,{1\over\Gamma(q) u}\,\times$& \cr
\+ & &\qquad\hfill $({u\over\upsilon})^{q+1}\, exp(-{u\over\upsilon})\,d\upsilon]\,dx$& \cr
\+ &=&$1-\sum_{m=0}^{+\infty}\int_0^t \int_{\R^+_*}{\theta^m\,u^q\,x^{p+m-1}\over m!\Gamma(p+m)\Gamma(q)}({1\over\upsilon})^{p+q+2m+1}exp(-{\theta+u+x\over\upsilon}) d\upsilon dx$& \cr
\+ &=&$1-\sum_{m=0}^{+\infty}\int_0^t{\Gamma(p+q+2m)\over m!\Gamma(p+m)\Gamma(q)}\,{\theta^m\,u^q\,x^{p+m-1}\over (\theta+u+x)^{p+q+2m}}\,dx$& \cr
\+ &=&$1-\sum_{m=0}^{+\infty}\int_0^{t\over \theta+u}{\Gamma(p+q+2m)\over m!\Gamma(p+m)\Gamma(q)}\,{\theta^m\,u^q\over (\theta+u)^{q+m}}{y^{p+m-1}\over (1+y)^{p+q+2m}}\,dy$& \cr
\+ &=&$1-\sum_{m=0}^{+\infty}{\Gamma(q+m)\over m!\Gamma(q)}\,{\theta^m\,u^q\over (\theta+u)^{q+m}}\,F(p+m,q+m,{t\over \theta+u})$& \cr
\dli $F(p+m,q+m,{t\over \theta+u})$ \'etant la valeur en ${t\over \theta+u}$ de la fonction de r\'epartition d'une loi $\beta(p+m,q+m)$ sur $\R^+$.
Pour $\theta=0$ on obtient une masse $Q^{(t,u)}([0])=1-F(p,q,{t\over u})$.
Dans le cas de l'analyse de variance, $p={k\over 2}$ et $q={l\over 2}$, cette masse est \'egale \`a la probabilit\'e qu'une loi de Fisher, de degr\'es de libert\'e $k$ et $l$, soit sup\'erieure \`a ${t/k\over u/l}$. C'est le seuil minimum de rejet du test de $H_0 : \theta=0$ contre $H_1 : \theta>0$. Nous avions d\'ej\`a obtenu ce type de r\'esultat pour le param\`etre ${\theta\over\sigma^2}$ dans l'exemple 3.3.
\dli Ce qui pr\'ec\`ede nous permet de d\'efinir une probabilit\'e inductive $Q^{(t,u)}$ sur $\Theta=\R^+$ puisque pour tout $\upsilon\in\Upsilon$ la proposition 4.2 s'applique (la fonction de r\'epartition $F_\upsilon(\theta,t)$ tend bien vers $0$ lorsque $\theta$ tend vers $+\infty$ car  $\Gamma(p+m,\upsilon,t)=\Gamma(p+m,1,{t\over\upsilon})$ tend vers $0$ lorsque $m$ tend vers $+\infty$). On peut ainsi traiter d'autres hypoth\`eses que les hypoth\`eses unilat\'erales.

Comme en 3.3 les expertises obtenues dans cet exemple compl\`etent l'infor\-mation donn\'ee par les r\'esultats des tests classiques, mais cette fois directement sur la valeur du param\`etre \'etudi\'e et non pas sur sa valeur exprim\'ee en \'ecart-type.  Le seuil minimum de rejet $\alpha_m({t\over u})$ du test de $H_0 : \theta=0$ contre $H_1 : \theta>0$ est vu comme le vote $Q^{(t,u)}([0])$ des experts en faveur de l'hypoth\`ese $H_0 : \theta=0$. 
Dans le cas du non rejet de $H_0$, c'est-\`a-dire lorsque $\alpha_m({t\over u})$ est sup\'erieur au seuil choisi, on peut pr\'eciser cette r\'eponse en regardant si $Q^{(t,u)}([0,\theta])$ se rapproche rapidement de $1$ lorsque $\theta$ cro\^{\i}t. Le cas du rejet de $H_0$ pose probl\`eme lorsque l'hypoth\`ese $H_0 : \theta=0$ est une id\'ealisation de l'hypoth\`ese r\'eelle \`a tester. Bien souvent l'utilisateur se demande si $\theta$ est petit et non pas si $\theta$ est nul. Une interpr\'etation trop rapide du rejet peut conduire \`a consid\'erer que $\theta$ est notable alors qu'il est n\'egligeable. Une analyse de la croissance de $Q^{(t,u)}([0,\theta])$ lorsque $\theta$ s'\'eloigne de $0$ permet d'\'eviter facilement ce pi\`ege. Il est cependant plus satisfaisant d'essayer de traduire l'hypoth\`ese ``$\theta$ est petit'' par $\theta\in[0,\theta_0]$ et de porter un jugement \`a partir de $Q^{(t,u)}([0,\theta_0])$.
\dli La prise en compte des votes $Q^{(t,u)}([0,\theta])$ autour de l'hypoth\`ese nulle $\theta=0$ donne une mani\`ere simple et parlante de compl\'eter l'infor\-mation donn\'ee par les tests classiques. Elle est moins complexe \`a analyser que la notion de puissance et de plus elle ne d\'epend pas de l'\'ecart-type inconnu $\sigma$. 

\ESendexample

\bigskip\bigskip
\ESsection{Conclusion}%
\medskip

Dans des mod\`eles \`a rapport de vraisemblance monotone sur $\R$ nous avons vu que la {\it p}-value associ\'ee \`a un test unilat\'eral est g\'en\'eralement consid\'er\'ee comme un indice de la confiance qu'on accorde \`a l'hypoth\`ese $H_0$. Avec notre approche cet indice s'interpr\`ete comme le r\'esultat d'un vote des experts en faveur de $H_0$, ce vote \'etant d\'efini sans information a priori.
La notion d'expertises compatibles pour l'ensemble des hypoth\`eses unilat\'erales nous a conduit \`a prolonger ces votes en une probabilit\'e inductive $Q^x$ sur l'espace des param\`etres, $x$ \'etant la r\'ealisation obtenue.
On retrouve ainsi la distribution fiduciaire de Fisher. 
La probabilit\'e $Q^x$ est souvent une loi a posteriori associ\'ee \`a une mesure a priori
qui peut \^etre consid\'er\'ee comme non informative dans le cadre bay\'esien.

Prendre en compte l'ensemble des hypoth\`eses unilat\'erales n'a de sens que si l'ordre sur
l'espace des param\`etres est totalement structurant pour l'interpr\'etation. Dans ce cas nous avons propos\'e, pour des raisons de coh\'e\-rence, de prendre comme indice de confiance d'une hypoth\`ese de la forme $\theta\in[\theta_1,\theta_2]$ la valeur sur cet intervalle de la probabilit\'e induite $Q^x$. C'est le prolongement aux hypoth\`eses bilat\'erales du vote des experts unilat\'eraux. Ce n'est pas la {\it p}-value du test bilat\'eral associ\'e \`a $H_0 : \theta\in[\theta_1,\theta_2]$. Cette mani\`ere de faire permet de prendre en compte les trois hypoth\`eses : $]-\infty,\theta_1[$, $[\theta_1,\theta_2]$, $]\theta_2,+\infty[$ , et d'avoir des interpr\'etations diff\'erentes pour chacune d'elles. Les tests bilat\'eraux ne sont pas construits pour cela. 
\dli Les probl\`emes soulev\'es par les {\it p}-values des deux tests bilat\'eraux d\'efinis \`a partir de l'hypoth\`ese $\theta\in[\theta_1,\theta_2]$ proviennent de l'\'equivalence d'interpr\'eta\-tion entre $\theta<\theta_1$ et $\theta>\theta_2$ que ces tests supposent. Dans l'exemple des lois $N(\theta,1)$ ou $N(\theta,\sigma^2)$, nous avons retrouv\'e leurs {\it p}-values comme votes d'experts en ne consid\'erant que les \'ev\'enements sym\'etriques par rapport \`a $c=(\theta_1+\theta_2)/2$. Ces {\it p}-values nous semblent recommandables uniquement dans le cas o\`u l'interpr\'etation est la m\^eme que $\theta$ soit inf\'erieur \`a $\theta_1$ ou sup\'erieur \`a $\theta_2$. La structure du probl\`eme de d\'ecision change. L'ordre sur $\Theta$ n'est plus alors totalement structurant pour l'interpr\'etation, c'est la distance \`a $c$ qui devient primordiale. Les hypoth\`eses l\'egitimes sont inva\-rian\-tes pour cette distance, elles sont sym\'etriques par rapport \`a $c$ et leurs {\it p}-values sont coh\'erentes. On peut g\'en\'eraliser cette mani\`ere de faire aux mod\`eles exponentiels r\'eels. La notion de coh\'erence des votes nous impose de traiter globalement un ensemble de probl\`emes de d\'ecision, sur un mod\`ele statistique donn\'e.
\dli Le choix entre les deux types de traitement des hypoth\`eses bilat\'erales d\'epend de la mani\`ere dont l'utilisateur structure son interpr\'etation des diff\'erentes valeurs du param\`etre. C'est la structuration par l'ordre sur $\Theta\subseteq\R$ qui est la plus courante. Dans ce cas les {\it p}-values des tests bilat\'eraux ne devraient pas \^etre utilis\'ees.

Au d\'epart la notion d'expertise n'a pas \'et\'e construite pour traiter le probl\`eme des {\it p}-values mais pour essayer de chercher une aide \`a la d\'ecision sous la forme d'un indice de confiance pour chacune des hypoth\`eses. Ce type d'inf\'erence a \'et\'e d\'efendu dans de nombreux travaux (Kroese, van der Meulen, Poortema et Schaafsma (1995)). Il est particuli\`erement adapt\'e aux applications qui ne r\'eclament pas un choix d\'efinitif entre les deux hypoth\`eses. Les solutions les plus courantes se font dans le cadre bay\'esien ou celui de l'estimation. Nous avons essay\'e de nous passer du choix souvent difficile d'une loi a priori ou d'une fonction de perte.
Le cas particulier des votes neutres, compatibles sur les hypoth\`eses unilat\'erales, donne des r\'esultats qui sont \`a rapprocher de ceux de l'analyse fiduciaire de Fisher. Utilis\'ee sans information a priori, la notion d'expertise s'inscrit dans le champ de questionnement ouvert par ce que l'on peut appeler le compromis de Fisher-Neyman-Wald.

\ESbiblio{Bibliographie}

\par%
{\sc Barra J. R.}
(1971), 
{\sl Notions fondamentales de statistique math\'ematique},
Dunod, Paris.

\par%
{\sc Berger, J. O.}
(1985), 
{\sl Statistical decision theory and Bayesian analysis} (second edition), Springer-Verlag, New York.

\par%
{\sc Berger, J. O.} et {\sc Sellke, T.}
(1987), 
Testing a point null hypothesis: the irreconcilability of {\sl p}-values and evidence, 
{\sl J. Amer. Statist. Assoc.} 
{\bf 82} 112--122.

\par%
{\sc Buehler, R. J.}
(1980), 
Fiducial inference. {\sl R. A. Fisher : an appreciation}, Springer-Verlag, New York, Lecture Notes in Statistics 109-118.

\par%
{\sc Casella, G.} et {\sc Berger, R. L.}
(1987), 
Reconciling Bayesian and frequentist evidence in the one-sided testing problem, 
{\sl J. Amer. Statist. Assoc.} 
{\bf 82} 106--111.

\par%
{\sc Dickey, J. M.}
(1968), 
Three multidimensional-integral identities with bayesian applications, 
{\sl Ann. Math. Statist.} 
{\bf 39} 1615--1628.

\par%
{\sc Fisher, R. A.}
(1935), 
The fiducial argument in statistical inference, 
{\sl Annals of Eugenics} 
{\bf 6} 391--398.

\par%
{\sc Gabriel, K. R.}
(1969), 
Simultaneous test procedures $-\!-$ some theory of multiple comparisons, 
{\sl Ann. Math. Statist.} 
{\bf 40} 224--250.

\par%
{\sc Hung, H. M. J.}, {\sc O'Neill, R. T.}, {\sc Bauer, P.} et {\sc K\"ohne, K.}  
(1997), 
The behavior of the {\sl p}-value when the alternative hypothesis is true, 
{\sl Biometrics} 
{\bf 53} 11--22.

\par%
{\sc Hwang, J. T.}, {\sc Casella, G.}, {\sc Robert, C.}, {\sc Wells, M. T.}, {\sc Farrell, R. H.} 
(1992), 
Estimation of accuracy in testing, 
{\sl Ann. Statist.} 
{\bf 20} 490--509.

\par%
{\sc Karlin, S.} 
(1955), 
Decision theory for P\'olya type distributions. Case of two actions, I, 
{\sl Proc. Third Berkeley Symposium on Math. Statist. and Prob.} (Univ. of Calif. Press) 
{\bf 1} 115--128.

\par%
{\sc Kroese, A. H.}, {\sc van der Meulen, E. A.}, {\sc Poortema, K.}, {\sc Schaafsma, W.} 
(1995), 
Distributional inference, 
{\sl Statistica Neerlandica} 
{\bf 49} 63--82.

\par%
{\sc Lehmann, E. L.} 
(1986),
{\sl Testing statistical hypotheses} (second edition), 
Wiley, New York.

\par%
{\sc van der Meulen, E. A.}, {\sc Schaafsma, W.} 
(1993), 
Assessing weights of evidence for discussing classical statistical hypotheses, 
{\sl Statistics \& Decisions} 
{\bf 11} 201--220.

\par%
{\sc Monfort, A.} 
(1982),
{\sl Cours de statistique math\'ematique}, 
Economica, Paris.

\par%
{\sc Morel, G.} 
(1997), 
Expertises : proc\'edures statistiques d'aide \`a la d\'ecision, 
{\sl Pr\'e-publi\-ca\-tion}, 182 p., 
LAST-Universit\'e de Tours. 

\par%
{\sc Morel, G.} 
(1998), 
Probabiliser l'espace des d\'ecisions, 
{\sl Pub. S\'em. 97 : Bru Huber Prum,} 
Paris V, 71--98. 

\par%
{\sc Robert, C.} 
(1992), 
{\sl L'analyse statistique bay\'esienne}, 
Economica, Paris.

\par%
{\sc Salom\'e, D.} 
(1998), 
{\sl Statistical inference via fiducial methods}, 
PhD thesis, Rijksuniversiteit Groningen.

\par%
{\sc Schaafsma, W.}, {\sc Tolboom, J.} et {\sc van der Meulen, B.} 
(1989), 
Discussing truth and falsity by computing a Q-value, in 
{\sl Statistical Data Analysis and Inference} (Y. Dodge, ed.), North-Holland, Amsterdam, 85--100. 

\par%
{\sc Scheff\'e H.}
(1970), 
{\sl The analysis of variance} (sixi\`eme \'edition), 
Wiley, New York.

\par%
{\sc Schervish, M. J.} 
(1996), 
{\sl P}-values : what they are and what they are not, 
{\sl Amer. Statist.} 
{\bf 50} 203--206.

\par%
{\sc Thompson, P.} 
(1996), 
Improving the admissibility screen : evaluating test statistics on the basis of p-values rather than power, 
{\sl Comm. Statist. Theory Methods} 
{\bf 25} 537--553.

\ESendbiblio 

\bye